\documentclass[a4paper]{article}
\usepackage{amsmath}
\usepackage{amssymb,latexsym}
\usepackage[dvips]{graphicx}
\usepackage{color}
\usepackage{cite}
\usepackage{amsthm}
\usepackage[affil-it]{authblk}

\usepackage{etoolbox}
\makeatletter
\patchcmd{\@maketitle}
  {\ifx\@empty\@dedicatory}
  {\ifx\@empty\@date \else {\vskip3ex \centering\footnotesize\@date\par\vskip1ex}\fi
   \ifx\@empty\@dedicatory}
  {}{}
\patchcmd{\@adminfootnotes}
  {\ifx\@empty\@date\else \@footnotetext{\@setdate}\fi}
  {}{}{}
 
\makeatother

\newtheorem{theorem}{Theorem} [section]
\newtheorem{lemma}{Lemma} [section]
\newtheorem{proposition}{Proposition} [section]

\newtheorem{definition}{Definition} [section]

\newtheorem{remark}{Remark}[section]

\let\ssection=\section\renewcommand{\section}{\setcounter{equation}{0}\ssection}

\begin{document}

\title{On the stability of the solitary waves to the rotation Benjamin-Ono equation.}
%\address{M. Darwich: Universit\'e Fran\c{c}ois rabelais de Tours, Laboratoire de Math\'ematiques
%et Physique Th\'eorique, UMR-CNRS 7350, Parc de Grandmont, 37200
%Tours, France} \email{Mohamad.Darwich@univ-tours.fr}
\date{}
%\predate{}
%\ofoot{\today}
\author{Mohamad Darwich.}
  %\thanks{Electronic address: \texttt{Mohamad.Darwich@lmpt.univ-tours.fr}.}
\affil{ Faculty of Sciences, Lebanese University Hadat-Lebanon.}
\maketitle
\begin{abstract}
In this paper, we study several aspects of solitary wave solutions of the rotation 
 Benjamin-Ono equation. By solving a minimization problem on the line, we construct a family of even travelling waves $\psi_{c,\gamma}$. We also study the strong convergence of this family and we establish the uniqueness of $\psi_{c,\gamma}$ for $\gamma$ small enough. Note that this improves the results in \cite{Esfahani} where the stability of the set of ground states is proven. 
\end{abstract}
%\author{ Mohamad Darwich.}

\maketitle

\section{Introduction}
In this paper we consider the rotation Benjamin-Ono (RBO) equation which can be written as:
\begin{equation}\label{RMBO}
(u_t - \mathcal{H}u_{xx} + (u^2)_x)_x + \gamma u=0, x \in \mathbb{R}, t>0
\end{equation}
where $\gamma >0$ is a real constant and $\mathcal{H}$ denotes the Hilbert transform defined by:
$$
\mathcal{H} u(t,x) = p.v.\frac{1}{\pi}\int_{\mathbb{R}} \frac{u(y,t)}{x-y}dy
$$
where p.v is the Cauchy principal value.
Equation \ref{RMBO} models the propagation of long internal waves in a deep rotation fluid see\cite{Galkin}, \cite{Linares},\cite{Redekopp} and \cite{Grimshaw} where u(t,x) represents the free surface of the liquid and the parameter $\gamma$ measures the effect of the rotation.% The parameter $\beta$ determines the type of dispersion.\\
\\The Cauchy problem of (\ref{RMBO}) is globally well-posed in $H^{\frac{1}{2}}(\mathbb{R})$ (see for instance \cite{Kenig}). Moreover, equation (\ref{RMBO}) has the following conserved quantities:
$$\displaystyle{\text{Energy:} ~~E(u(t))= \int_{\mathbb{R}}\left(\frac{1}{2}(D_x^{\frac{1}{2}}u)^2 + \frac{\gamma}{2}(\partial_x^{-1}u)^2 - \frac{1}{3}u^3\right)dx}=E(u(0)),$$
$$\displaystyle{\text{Mass}:~V(u(t))=\int_{\mathbb{R}}u^2dx=V(u(0))}.$$

We are interested in solitary waves of (RBO); i.e, the solutions to equation (\ref{RMBO}) of the form $u(t,x) = \psi(x-ct)$ traveling with the speed $c \in \mathbb{R^{+}}$. By substituting $u$ by $\psi$ in (\ref{RMBO}), integrating on $\mathbb{R}$, we obtain:
\begin{equation}\label{solitaryequationRBO}
\mathcal{H}\partial_x\psi + c \psi -\gamma \partial_x^{-2}\psi= \psi^2
\end{equation}
In \cite{Esfahani}, Esfahani and Levandosky have studied the solitary waves and the orbital stability of the (RBO) equation, and up to our knowledge no uniqueness result (up to symmetries) is known even for ground state solutions.\\
However, if we take $\gamma =0$ in (\ref{RMBO}) and integrate with respect to $x$ we obtain the Benjamin-Ono (BO) equation:
\begin{equation}\label{MBO}
u_t - \mathcal{H}u_{xx} + (u^2)_x = 0,
\end{equation}
and thus for $\gamma = 0$ in (\ref{solitaryequationRBO}) we recover the equation of solitons of (BO):

\begin{equation}\label{solitaryequationBO}
\mathcal{H}\partial_x\psi + c \psi = \psi^2.
\end{equation}
It is known that for every $c>0$ there is a unique (up to translations) solution of (\ref{solitaryequationBO})
which is:
\begin{equation}\label{solitonBO}
Q_{c}(x)=\frac{4c}{1+c^2x^2},
\end{equation}
(see Benjamin \cite{Benjamin} and Amick and Toland \cite{Amick} for the uniqueness statement).
This solution is stable in a sense which will be defined later (see Bennet et al.\cite{Bennett} and Weinstein \cite{Weinstein}).\\
In this paper, inspired from \cite{Molinet}, we minimize for all $c\in \mathbb{R^+}$ and $\gamma >0$ in the even functions the functional:
\begin{equation}
\displaystyle{I(u) = I(u,c,\gamma) = \int_{\mathbb{R}}\left( \frac{1}{2}(\partial_x^{\frac{1}{2}}u)^2 + \frac{c}{2}u^2  +\frac{\gamma}{2}(\partial_x^{-1}u)^2\right)dx}
\end{equation}
under the constraint:
\begin{equation}\label{constraint}
K(u) =\frac{1}{3}\displaystyle{\int u^3dx=K(Q_{c})},
\end{equation}
 we also prove the uniqueness of the associated even ground states and establish that 
the solitary waves of the (BO) equation are the limits in $H^{\frac{1}{2}}$ of the even ground states of the (RBO) equation. Finally, we prove that the even ground states of (RBO) are stable in $H^{\frac{1}{2}}(\mathbb{R})$.
\begin{remark}
We chose this constraint in (\ref{constraint}) to be sure that the solitary waves of the (RBO) equation will converge to the solitons $Q_{c}$.\end{remark}
Let us first recall the definition of orbital stability and define what we will call ground state solutions to (\ref{solitaryequationRBO}):
\begin{definition}
Let $\psi \in H^{\frac{1}{2}}(\mathbb{R})$ be a solution of (\ref{solitaryequationRBO}). We say that $\psi$ is orbitally stable in $H^{\frac{1}{2}}(\mathbb{R})$, if for all $\epsilon >0$, there exists $\delta>0$ such that for all initial data $u_0 \in H^{\frac{1}{2}}(\mathbb{R})$, satisfying 
$$
||u_0-\psi||_{H^{\frac{1}{2}}} \leq \delta,
$$
the solution $u \in C(\mathbb{R_{+}},H^{\frac{1}{2}})$ of (\ref{RMBO}) emanating from $u_0$ satisfies:$$
\displaystyle{\sup_{t\in\mathbb{R_{+}}}\inf_{z\in\mathbb{R}}||u(t,.+z)-\psi||_{H^{\frac{1}{2}}} \leq \epsilon.}$$
\end{definition}
Let $X$ be the space defined by: 
$$
X= \left\{f \in H^{\frac{1}{2}}(\mathbb{R}); \partial_x^{-1}f \in L^2(\mathbb{R})~~\text{and}~~f(-x)=f(x)~\forall x \in \mathbb{R}\right\},
$$
with the norm:
$$
||f||_{X} = ||f||_{H^{\frac{1}{2}}(\mathbb{R})} + ||\partial_x^{-1}f||_{ L^2(\mathbb{R})}.
$$
This space endowed with its metric is an Hilbert space.
\begin{definition}
We say that a solution to (\ref{solitaryequationRBO}) is an even ground state solution to (\ref{solitaryequationRBO}), if it is also a solution to the constraint minimization problem
$$
\displaystyle{M_{\lambda} = \inf\{I(u), u \in X, K(u) = \lambda\}}
$$
for some $\lambda >0.$
\end{definition}
This paper will be organized as following: in Section 2 we will prove the existence of even solitary vaves $\psi_{c,\gamma}$ of (\ref{RMBO}), in Section 3 we study the strong convergence in $H^{\frac{1}{2}}(\mathbb{R})$ of the family $\{\psi_{c,\gamma}, 0<\gamma <<1\}$ to the explicit solitons $Q_{c}$ of (\ref{solitonBO}) as $\gamma$ goes to 0 and in the last section we establish the uniqueness of $\psi_{c,\gamma}$ for $\gamma$ small enough and the orbital stability in $H^{\frac{1}{2}}$.\\
Let us now give our main results:
\begin{theorem}\label{theorem}
Let $\gamma_0 >0$ fixed.  There exists $\delta>0$ such that for all $c > \frac{\gamma_0}{\delta}$, there exists a unique even ground state solution $\psi_{c,\gamma_0}$ to equation (\ref{solitaryequationRBO}). These solutions form a $H^{\frac{1}{2}}$-continuous curve of even solitary waves of (\ref{RMBO}) which are all stable.
\end{theorem}
%Acknowledgments : The  author thank the L.M.P.T. for his kind  hospitality during the development of this work. Moreover, he would like to thank Professor Luc Molinet for valuable remarks and comments.
\section{Existence of solitary waves of (RBO) equation}
In this section, we prove the existence of even solitary waves of (RBO) for $\gamma>0$, by solving a minimization problem. In this step we argue as in \cite{Levandosky} and \cite{Esfahani} and we need the following lemma: % We start by proving that the functional $I$ is coercive in the space $X$:% to prove this we need the following lemma:
%\begin{lemma}
%If $\psi$ is a solution of (\ref{solitaryequationRBO}), then:$$
%\displaystyle{2\int_{\mathbb{R}}(\partial_x^{\frac{1}{2}}\psi)^2dx - c \int_{\mathbb{R}}\psi^2 + 7\gamma \int_{\mathbb{R}}(\partial_x^{-1}\psi)^2dx=0}.$$
%\end{lemma}
%\proof See Lemma 3.1 in \cite{Esfahani}.
\begin{lemma}\label{coercive}
Let $\gamma>0$ and $c>0$, then the  functional $I$ is coercive in the space $X$.
\end{lemma}
\proof We have that
\begin{equation}\label{eq1}
 \min {(\frac{1}{2},\frac{c}{2},\frac{\gamma}{2})}||u||^2_X\leq I(u)
\end{equation}
then $I$ is coercive on $X$.
%\geq \int_{\mathbb{R}}(\partial_x^{\frac{1}{2}}u)^2 + \gamma(\partial_x^{-1}u)^2dx}
\\
%\begin{equation}\label{eq2}
%\max{(c,\gamma) \|u\|^2_X \geq I(u) \geq \int_{\mathbb{R}}c_1(\partial_x^{\frac{1}{2}}u)^2 + %c_2\gamma(\partial_x^{-1}u)^2dx}
%\end{equation}
%where $c_1 = 1 -\sqrt{\frac{1}{2} + \frac{2c^3}{27\gamma}}$ and $c_2 = \frac{27 - 4c^3}{27\gamma + 4c^3}.$
%Note that (\ref{eq1}) and (\ref{eq2}) and the following inequality:
%\begin{equation}\label{AB}
%\|u\|_{L^2(\mathbb{R})}^{2} \leq A \|\partial_x^{\frac{1}{2}}u\|_{L^2(\mathbb{R})}^{2} + B \|\partial_x^{-1}u\|_{L^2(\mathbb{R})}^{2}
%\end{equation}
%where $A>0$, $B>\frac{4}{27A^2}$. %we obtain the coercivity condition $I(u) \sim ||u||^2_X$.\\
%Moreover,  
%\begin{align}I(u) &=\int_{\mathbb{R}} \frac{1}{2}(\partial_x^{\frac{1}{2}}u)^2 + \frac{c}{2}\int u^2 + \frac{\gamma}{2}(\partial_x^{-1}u)^2 -\gamma\left(\int (\partial_x^{-1}u)^2dx \right)\nonumber
%\end{align}
 %and from the following inequality:
%\begin{equation}\label{AB}
%\|u\|_{L^2(\mathbb{R})}^{2} \leq A \|\partial_x^{\frac{1}{2}}u\|_{L^2(\mathbb{R})}^{2} + B \|\partial_x^{-1}u\|_{L^2(\mathbb{R})}^{2},~~ \text{where}~~ A>0, B>\frac{4}{27A^2}.
%\end{equation}
%we infer that for $\gamma>0$ small, $I$ is coercive in $X$.\\

For $\lambda >0$ given, we define $M_{\lambda}$ as $\displaystyle{M_{\lambda}=\inf\left\{I(u), u \in X, K(u) = \lambda\right\},}$ $M_{\lambda}$ satisfies the following property:
\begin{lemma}\label{Mpositive}
Let $\lambda >0$ given, then $M_{\lambda} >0.$
\end{lemma}
\proof 

By Sobolev embedding and by interpolation we obtain that:
\begin{equation}\label{eq3}
\|u\|^3_{L^3(\mathbb {R})} \leq C \|u\|^3_{H^{1/6}(\mathbb {R})} \leq C \|u\|^{5/3}_{H^{\frac{1}{2}}(\mathbb {R})}\|u\|^{4/3}_{H^{-1/4}(\mathbb {R})}.
\end{equation}
Now by Cauchy-Schwarz we have that:
\begin{equation}\label{eq4}
\|u\|_{H^{-1/4}(\mathbb {R})} \leq C \|u\|^{\frac{1}{2}}_{H^{\frac{1}{2}}(\mathbb {R})}||\partial_x^{-1}u||_{L^2(\mathbb {R})}^{\frac{1}{2}}.
\end{equation}
From  (\ref{eq3}) and (\ref{eq4}) we obtain that:

\begin{equation}\label{eq5}
\|u\|^3_{L^3(\mathbb {R})} \leq C \|u\|_{H^{\frac{1}{2}}(\mathbb{R})}^{7/3}\|\partial_x^{-1}u\|_{L^2(\mathbb{R})}^{\frac{2}{3}}.
\end{equation}
Thus we can deduce that:
\begin{align}
\lambda = K(u) \leq C \int |u|^3dx &\leq C\|u\|_{H^{\frac{1}{2}}(\mathbb{R})}^{7/3}\|\partial_x^{-1}u\|_{L^2(\mathbb{R})}^{\frac{2}{3}}\nonumber\\
&\leq C \left( ||u||^2_{L^2(\mathbb{R})} + ||\partial_x^{\frac{1}{2}}u||^2_{L^2(\mathbb{R})} + ||\partial_x^{-1} u ||^2_{L^2(\mathbb{R})}\right)^{3/2}\nonumber
\end{align}
this gives that, $\lambda^{\frac{2}{3}} \leq C I(u)$, for all $ u \in X$ then $M_{\lambda} \geq \left(\frac{\lambda^{\frac{2}{3}}}{C}\right) >0$.\\

Now we will prove the existence of even ground states of the (RBO) equation, more precisely we have the following proposition:
\begin{proposition}\label{minimizing}
Let $\gamma >0$ and $c>0$ and $\{\psi_n\}_n$ be a minimizing sequence for some $\lambda >0$. Then there exist a subsequence and scalars $y_n \in \mathbb{R}$ and $\psi \in X$ such that $\psi_n(.+y_n) \longrightarrow \psi $ in $X$. Moreover, the function $\psi$ achieves the minimum $I(\psi) = M_{\lambda}$ subject to the constraint $K(\psi) = \lambda$.

\end{proposition}
\proof
We argue similary as Levandosky \cite{Levandosky}. %By proposition \ref{coercive}, we have that $I$ is coercive in $X$, let us prove now that $K$ is locally Lipschitz in $L^3(\mathbb{R})$. Indeed, for all $(u,v) \in L^3(\mathbb{R})\times L^3(\mathbb{R})$ applying the Holder inequality, we have:
%\begin{align}\label{Lip}
%\left|K(u)-K(v)\right|& = \left|\int_{\mathbb{R}}(u^3-v^3)dx\right|\nonumber\\
%&\leq ||u-v||_{L^3(\mathbb{R})}\sum_{k=0}^{2}|||u|^{2-k}|v|^{k}||_{L^{3/2}(\mathbb{R})},
%\end{align}
%and  by Young inequality, we have:
%\begin{equation}\label{Lip1}
%|||u|^{2-k}|v|^{k}||^{3/2}_{L^{3/2}(\mathbb{R})}\leq\frac{2-k}{2}||u||^{3}_{L^3(\mathbb{R})} + \frac{k}{2}||v||_{L^3(\mathbb{R})}^3.
%\end{equation}
%Then combining (\ref{Lip}) and (\ref{Lip1}), we obtain 
%\begin{equation}
%|K(u)-K(v)| \leq C(||u||_{L^3(\mathbb{R})},||v||_{L^3(\mathbb{R})})||u-v||_{L^3(\mathbb{R})}.
%\end{equation}
We study the minimization of $I$ in $X$ subject to the constraint $K(,)=K(Q_{c})$. Note that we can easily see that $M_{\lambda} = \lambda^{\frac{2}{3}}M_1$ so that the strict subadditivity condition:
$$
M_{\alpha} + M_{\lambda - \alpha} > M_{\lambda},
$$
holds for any $\alpha \in (0,\lambda)$.
Let us solve the constraint minimization problem. We take $(\psi_k)_k \subset X$ a minimizing sequence of the problem, i.e, for all $k \in \mathbb{N^*}$, we have:
\begin{equation}\label{convenergie}
K(\psi_k) = K(Q_{c}) = \lambda ~\text{and} ~~\lim_{k \longrightarrow +\infty}|I(\psi_k)-M_{\lambda}|=0,
\end{equation}
From the convergence (\ref{convenergie}) and the coercivity of $I$ we can deduce that the sequence $(\psi_k)_k$ is bounded in $X$. Since $X$ is reflexive, there exists a sub-sequence $(\psi_j)_{j} \subset X$ and a function $\psi \in X$ such that $\psi_j \longrightarrow \psi$ weakly in $X$. Let 
$$\rho_n = |\partial_x^{\frac{1}{2}}\psi_n|^2 + |\partial_x^{-1}\psi_n|^2,$$\\
from (\ref{convenergie}), one can see that $(\rho_n)_n$ is bounded in $L^{1}(\mathbb{R})$. After extracting a subsequence, we may assume that $\displaystyle{\lim_{k\longrightarrow +\infty}\int \rho_n = L < +\infty.}$ By normalizing, we may assume further that $\displaystyle{\int \rho_n = L}$ for all $n\in \mathbb{N^*}$. Then by the Concentration-Compactness Lemma \cite{Levandosky}, there are three possibilities: vanishing, dichotomy or compactness conditions.

 In the same way as in\cite{Levandosky} and \cite{Levandosky1}, it follows from the coercivity of $I$, inequality (\ref{eq5}) and the subadditivity condition that both vanishing and dichotomy may be ruled out, and therefore the sequence $\rho_n$ is compact that is: there exists $y_k \in \mathbb{R}$ such that for any $\epsilon >0$ there exists $R(\epsilon)$ such that
 $$
 \int_{B(y_k,R(\epsilon))}\rho_k dx \geq \int_{\mathbb{R}} \rho_k dx - \epsilon, \forall k\in \mathbb{N}.
 $$ Now let $\phi_n(x) = \psi_n(x+y_n)$. Since $\phi_k$ is bounded in $X$, a subsequence $\phi_k$ converges weakly to some $\psi \in X$, and by the weak lower semicontinuity of I over $X$,  we have :
 $$I(\psi) \leq \lim_{k\longrightarrow +\infty} I(\phi_k) =M_{\lambda}.$$
 Furthemore, weak convergence in X, compactness of $\rho_n$ and the inequality (\ref{eq5}) imply strong covergence in $L^3$. Therefore
 $$
 K(\psi) = \lim_{k\longrightarrow +\infty} K(\phi_k) = \lambda,
 $$
so $I(\psi) \geq M_{\lambda}$. Together with the inequality above, this implie $I(\psi) = M_{\lambda}$, so $\psi$ is a minimizer of $I$ subject to the constraint $K(\phi) = \lambda$. Finally, since $I$ is equivalent to the norm on $X$, $\phi_k$ converges weakly to $\psi$ and $I(\phi_k) \longrightarrow I(\psi)$, it follows that $\phi_k$ converges strongly to $\psi$ in $X$.\\

Now let us define the following quantity:
$$
\displaystyle{m=m(c,\gamma)= \inf\left\{\frac{I(u)}{K(u)^{\frac{2}{3}}}, u \in X, K(u)>0 \right\}}.
$$
\begin{remark}\label{remark1} In \cite{Esfahani}, the authors prove that the minimum for  m are exactly the ground states of (\ref{solitaryequationRBO}) and remark that if we multiply (\ref{solitaryequationRBO}) by $\psi$ and integrate, we obtain that $I(\psi,c,\gamma) = K(\psi)$, thus we can deduce that if $\psi$ is a ground state of (\ref{solitaryequationRBO}) then $K(\psi)=I(\psi,c,\gamma)=(m(c,\gamma))^{3}$.
\end{remark}
Analogous to Proposition \ref{minimizing}, we have the following one:
\begin{proposition}\label{minimizing1}
Let $\psi_n \in H^{\frac{1}{2}}$  be a given sequence such that $\displaystyle{\lim_{n \longrightarrow +\infty} I(\psi_n,c,0)}$\\
$\displaystyle{= \lim_{n \longrightarrow +\infty} K(\psi_n)= (m(c,0))^{3}}$ then there exists a subsequence
 renamed $\psi_n$, scalars $y_n \in \mathbb{R}$ such that $ \psi_n(.+y_n) \longrightarrow Q_{c}$  in $H^{\frac{1}{2}}.$
\end{proposition}
\proof We will proceed in the same way as in Proposition \ref{minimizing}. %We study the minimization of $I$ in $X$ subject to the constraint $K(,)=K(Q_{c})$. Note that we can easily see that $M_{\lambda} = \lambda^{\frac{2}{3}}M_1$ so that the strict subadditivity condition:
%$$
%M_{\alpha} + M_{\lambda - \alpha} > M_{\lambda},
%$$
%holds for any $\alpha \in (0,\lambda)$.
By hypothesis $ I(\psi_n,c,0)$ is convergent, but $ ||\psi_n||^2_{H^{\frac{1}{2}}}\lesssim I(\psi_n,c,0)$ then $\psi_n$ is bounded in $H^{\frac{1}{2}}(\mathbb{R})$ and that there exists a subsequence (renamed $\psi_n$) and $\psi \in H^{\frac{1}{2}}$ such that $\psi_n \longrightarrow \psi$ weakly as $ n\longrightarrow +\infty$. 
Let 
$$\xi_n = |\partial_x^{\frac{1}{2}}\psi_n|^2 + |\psi_n|^2,$$ $(\xi_n)_n$ is bounded in $L^{1}(\mathbb{R})$, then after extracting a subsequence, we may assume that $\displaystyle{\lim_{k\longrightarrow +\infty}\int \xi_n = L < +\infty.}$ By normalizing, we may assume further that $\displaystyle{\int \xi_n = L}$ for all $n\in \mathbb{N^*}$. By the Concentration-Compactness Lemma and as in Proposition \ref{minimizing} the sequence $\xi_n$ is compact. Now let $v_n(x) = \psi_n(x+y_n)$, since $v_k$ is bounded in $H^{\frac{1}{2}}$, a subsequence $v_k$ converges weakly to some $v \in X$, and by the weak lower semicontinuity of I over $H^{\frac{1}{2}}$,  we have :
 $$I(v) \leq \lim_{k\longrightarrow +\infty} I(v_k) = m(c,0)^3.$$
 Now, weak convergence in $H^{\frac{1}{2}}$ , compactness of $\xi_n$ and the following inequality:
 \begin{equation}\label{L3L2} 
 ||u||_{L^3} \lesssim ||u||^{\frac{2}{3}}_{L^2}||u||^{\frac{1}{3}}_{H^{\frac{1}{2}}}
 \end{equation}
  imply strong covergence in $L^3$. Therefore
 $$
 K(v) = \lim_{k\longrightarrow +\infty} K(v_k) = m(c,0)^3.
 $$
The function $v$ satisfies $I(v)=K(v)=m(c,0)^3$ ($m(c,0)>0$ since $M_{\lambda} = m\lambda^{\frac{2}{3}}$), this ensures that $v=Q_c$ since the solitary waves of (BO) are unique. Finally, weak convergence and convergence in norm prove that $\psi_n(.+y_n) \longrightarrow Q_c$ in $H^{\frac{1}{2}}$, this ends the proof.%Together with the inequality above, this implie $I(\psi) = M_{\lambda}$, so $\psi$ is a minimizer of $I$ subject to the constraint $K(\phi) = \lambda$. Finally, since $I$ is equivalent to the norm on $X$, $\phi_k$ converges weakly to $\psi$ and $I(\phi_k) \longrightarrow I(\psi)$, it follows that $\phi_k$ converges strongly to $\psi$ in $X$.\\

\section{Convergence in $H^{\frac{1}{2}}(\mathbb{R})$}
In this section, we prove the strong convergence of the family of solitary waves to (\ref{solitaryequationRBO}) in $H^{\frac{1}{2}}(\mathbb{R})$ to the explicit solitons $Q_{c}$ defined in (\ref{solitonBO}).\\

More precisely, we have the following theorem:
\begin{theorem}\label{convergence}
Let $c>0$ and consider any sequence $\gamma_n \longrightarrow 0^+$. Denote by $\psi_n$ the associated ground states of (\ref{solitaryequationRBO}) constructed in
Proposition \ref{minimizing}. Then there exists a subsequence (renamed $\psi_n$) and translations $y_n$ so that
$$
\psi_n \longrightarrow Q_{c} ~\text{in}~ H^{\frac{1}{2}}
$$
as $\gamma_n \longrightarrow 0^+$. That is , the solitary waves of the GBO equation are the limits in $H^{\frac{1}{2}}$ of solitary waves of the RGBO equation.
\end{theorem}
Before proving Theorem \ref{convergence}, we need the following lemmas:
\begin{lemma}\label{density}
The space $X$ is dense in $H^{\frac{1}{2}}(\mathbb{R})$.
\end{lemma}
\proof Let $u \in H^{\frac{1}{2}}(\mathbb{R})$ and $n \in N^*$. We define $u_n$ as $\hat{u_n}(\xi)=\hat{u}(\xi)1_{|\xi| > \frac{1}{n}}(\xi).$ By Parseval's identity we have that:
$$
||\partial_x^{-1}u_n||_{L^2(\mathbb{R})}^2 = ||\xi^{-1}\hat{u_n}||_{L^2(\mathbb{R})}^2 = \int_{|\xi| > \frac{1}{n}}\xi^{-2}|\hat{u}(\xi)|^2 d\xi < n^2 ||u||^2_{L^2(\mathbb{R})} < + \infty.
$$
Since $||u_n||_{L^2(\mathbb{R})} \leq ||u||_{L^2(\mathbb{R})} < +\infty$ and since $||\partial_x^{\frac{1}{2}}u_n||_{L^2(\mathbb{R})} \leq ||\partial_x^{\frac{1}{2}}u||_{L^2(\mathbb{R})}< +\infty$, it follows that $u_n \in X$.
Now:
$$
||u_n-u||_{H^{\frac{1}{2}}(\mathbb{R})}^2 = \int_{|\xi| < \frac{1}{n}} (1+|\xi|)|\hat{u}(\xi)|^2d\xi \leq ||u||_{H^{\frac{1}{2}}(\mathbb{R})}^2 < +\infty
$$
then for $n$ large enough we obtain that:
$$
||u_n-u||_{H^{\frac{1}{2}}(\mathbb{R})}^2 = \int_{|\xi| < \frac{1}{n}} (1+|\xi|)|\hat{u}(\xi)|^2d\xi < \epsilon
$$
this ends the proof.

\begin{lemma}\label{m}
The function m is strictly increasing in $\gamma$.
\end{lemma}
\proof
Fix $0<c$ and let $ \gamma_1 < \gamma_2$.
%Let $\phi_{c_1}$ and $\phi_{c_2}$ be ground states with $c=c_1$ and 
\begin{align}
m(c,\gamma_1) &\leq \frac{I(\psi_{c,\gamma_2},c,\gamma_1)}{K(\psi_{c,\gamma_2})^{\frac{2}{3}}}\nonumber\\
&=\frac{I(\psi_{c,\gamma_2},c,\gamma_2) + (\gamma_1-\gamma_2)\int (\partial_x^{-1} \psi_{c,\gamma_2})^2}{{K(\psi_{c,\gamma_2})^{\frac{2}{3}}}}\nonumber\\
&=\frac{I(\psi_{c,\gamma_2},c,\gamma_2)}{K(\psi_{c,\gamma_2})^{\frac{2}{3}}} + \frac{(\gamma_1-\gamma_2)\int (\partial_x^{-1} \psi_{c,\gamma_2})^2}{K(\psi_{c,\gamma_2})^{\frac{2}{3}}}\nonumber\\
&=m(c,\gamma_2) +\frac{(\gamma_1-\gamma_2)\int (\partial_x^{-1} \psi_{c,\gamma_2})^2}{K(\psi_{c,\gamma_2})^{\frac{2}{3}}}\nonumber \\
&<m(c,\gamma_2).\nonumber
\end{align}
\begin{lemma}\label{continuityzero}
%Fix $0<c<c^*$. Then 
%$$
%\displaystyle{\lim_{\gamma \longrightarrow 0} m(c,\gamma) = m(c,0)}.
%$$
The function $m$ is continuous in  $\gamma = 0$.
\end{lemma}

%Let $0< \gamma$.
%\begin{align}
%m(c,\gamma) &\leq \frac{I(\psi_{c,0},c,\gamma)}{K(\psi_{c,0})^{\frac{2}{3}}}\nonumber\\
%&=\frac{I(\psi_{c,0},c,0) + \gamma\int (\partial_x^{-1} \psi_{c,0})^2}{{K(\psi_{c,0})^{\frac{2}{3}}}}\nonumber\\
%&=\frac{I(\psi_{c,0},c,0)}{K(\psi_{c,0})^{\frac{2}{3}}} + \frac{\gamma\int (\partial_x^{-1} \psi_{c,0})^2}{K(\psi_{c,0})^{\frac{2}{3}}}\nonumber\\
%&=m(c,0) +\frac{\gamma\int (\partial_x^{-1} \psi_{c,0})^2}{K(\psi_{c,0})^{\frac{2}{3}}}\nonumber \\
%\end{align}
%Now using the coercivity of $I$, we obtain that:
%$$
%0\leq m(c,\gamma)-m(c,0) \leq \gamma \min {(\frac{1}{2},\frac{c}{2},\frac{\gamma}{2})} m(c,0),
%$$
%so $m$ is continuous in $\gamma$.
\proof Since $m$ is strictly increasing in $\gamma$, it suffices to show that $m(c,\gamma_k) \longrightarrow m(c,0)$ for some sequence $\gamma_k \longrightarrow 0$. By Theorem \ref{density} we may choose a sequence $\psi_k$ in $X$ such that $||\psi_k-Q_{c}||_{H^{\frac{1}{2}}} < 1/k$ and let us define:

$$
\displaystyle{\gamma_k = \min\left(\frac{1}{k},\frac{1}{k\int|\partial_x^{-1}\psi_k|^2dx}\right).}
$$
Then 
\begin{align}
m(c,\gamma_k) &\leq \frac{I(\psi_k,c,\gamma_k)}{K(\psi_k)^{\frac{2}{3}}}\nonumber\\
&= \frac{I(\psi_k,c,0) + \gamma_k\int|\partial_x^{-1}\psi_k|^2dx}{K(\psi_k)^{\frac{2}{3}}}\nonumber\\
&\leq \frac{I(\psi_k,c,0) + 1/k}{K(\psi_k)^{\frac{2}{3}}}.\nonumber
\end{align}
Since $I(.,c,0)$ and $K$ are both continuous on $H^{\frac{1}{2}}$, we therefore have:
$$
\displaystyle{\limsup_{k\longrightarrow +\infty}m(c,\gamma_k) \leq \frac{I(Q_{c},c,0)}{K(Q_{c})^{\frac{2}{3}}}= m(c,0).}
$$
On the other hand, let $\psi_n$ be the ground states of RGBO for $\gamma=\gamma_n$ for every $n \in N$, so
\begin{align}
m(c,0) &\leq \frac{I(\psi_n,c,0) }{K(\psi_n)^{\frac{2}{3}}}\nonumber\\
&=\frac{I(\psi_n,c,\gamma_n) -\gamma_n\int|\partial_x^{-1}\psi_n|^2dx}{K(\psi_n)^{\frac{2}{3}}}\nonumber\\
&< \frac{I(\psi_n,c,\gamma_n)}{K(\psi_n)^{\frac{2}{3}}} = m(c,\gamma_n).\nonumber
\end{align}
Thus 
$$
\displaystyle{\liminf_{n\longrightarrow+\infty} m(c,\gamma_n) \geq m(c,0),}
$$
and the lemma follows.\\

Let us prove theorem \ref{convergence}: Let $\psi_n$ be the solitary waves of RGBO for $\gamma=\gamma_n$ for every $n \in N$.
By remark (\ref{remark1}) and continuity of m at $\gamma=0$ we have:$$
\displaystyle{\lim_{n \longrightarrow +\infty} I(\psi_n,c,\gamma_n)=\lim_{n \longrightarrow +\infty} K(\psi_n) =\lim_{n\longrightarrow +\infty}m(c,\gamma_n)^{3} = m(c,0)^{3}}
$$
and 
\begin{align}
\lim_{n\longrightarrow +\infty} I(\psi_n,c,0) &= \lim_{n\longrightarrow +\infty}\left( I(\psi_n,c,\gamma_n) -\gamma_n\int(\partial_x^{-1}\psi_n)^2dx\right)\nonumber\\
&\leq \lim_{n\longrightarrow +\infty} I(\psi_n,c,\gamma_n)\nonumber\\
&=\lim_{n\longrightarrow +\infty} m(c,\gamma_n)^3\nonumber\\
&=m(c,0)^3.\nonumber
\end{align}

Now the result follows from Proposition \ref{minimizing1}.
\section{Uniqueness of even ground states and the orbital stability result}

In this section, we fix $c=1$, and we will prove that for $\gamma$ small enough, the even ground state solution to (\ref{solitaryequationRBO}), with $c=1$, is unique. More precisely we have the following lemma:
\begin{lemma}
Let $\gamma>0$, and $c=1$. There exists $\delta >0$ such that for all $\gamma \in]0,\delta[$, there exists a unique even ground state solution $\psi_{1,\gamma}$ to equation (\ref{solitaryequationRBO}). Moreover, the map $\gamma \longmapsto \psi_{1,\gamma}$ is continuous from $]0,\delta[$ into $H^{\frac{1}{2}}(\mathbb{R})$.
\end{lemma}
\begin{remark}\label{dilatation}
The previous lemma gives the uniqueness result for $\gamma>0$ small and $c=1$. But remark that, $u$ satisfies (\ref{solitaryequationRBO}) with $\gamma>0$ and $c=1$ if and only if $v = \frac{\gamma}{\gamma_0} u(\frac{\gamma}{\gamma_0} x)$ satisfies (\ref{solitaryequationRBO}) with $\gamma_0>0$ and $c=\frac{\gamma_0}{\gamma}$. Then the uniqueness result holds for $\gamma_0>0$ and $c >0$ large enough.
\end{remark}
\proof Let us start by proving that for any $\epsilon >0$ there exists $\gamma_{\epsilon}>0 $ such that for any $\ 0<\gamma < \gamma_{\epsilon}$, any even ground states solution $\psi=\psi_{1,\gamma}$ to (\ref{solitaryequationRBO}) satisfies:
\begin{equation}\label{smallness}
||\psi - Q_{1}||_{H^{\frac{1}{2}}} < \epsilon.
\end{equation}
For this, we will proceed by contradiction, assume that there exist $\epsilon_0 >0$, a sequence of positive reals numbers $(\gamma_n)_n$ such that $\gamma_n  \longrightarrow 0$ as $n \longrightarrow 0$ and a sequence $(\psi_n)_n$ of ground states to (\ref{solitaryequationRBO}) such that:
\begin{equation}\label{contresmallness}
||\psi_n - Q_{1}||_{H^{\frac{1}{2}}} \geq \epsilon_0, \forall n \in \mathbb{N}.
\end{equation}
From equation (\ref{solitaryequationRBO}) we can write that $K(\psi_n) = I(\psi_n)>0$. Let
$$
w_n = (\frac{K(Q_{1})}{K(\psi_n)})^{1/3}\psi_n
$$
then $K(w_n) = K(Q_{1})$. Moreover, $(w_n)_{n}$ is a sequence of solutions to $M_{\lambda}$ and as in the proof of the convergence in $H^{\frac{1}{2}}$ in section 2, it follows that $w_n \longrightarrow Q_{1}$ in $H^{\frac{1}{2}}$. Therefore $$
\psi_n = \left(\frac{I(w_n)}{K(Q_{1})}\right)w_n \longrightarrow Q_{1},~\text{in}~~ H^{\frac{1}{2}}
$$
wich contradicts (\ref{smallness}).\\
\begin{remark}
Note that by (\ref{solitonBO}), we have that $||Q_{c}||_{L^2}^2 = 16 c D$, where $D>0$ is a constante does not depend on $c$. Then $\frac{d}{dc}||Q_{c}||_{L^2}^2 >0$ for any $c>0$.
\end{remark}
Now let us prove the uniqueness result. Let $\psi_{1,\gamma}$ and $\tilde{\psi}_{1,\gamma}$ be two ground states of (\ref{solitaryequationRBO}). Following the idea of Kenig et al \cite{Kenig1} (Proposition 3), we will prove, arguing by contradiction, that $\tilde{w} = \psi_{1,\gamma} - \tilde{\psi}_{1,\gamma}=0$ as soon as $\gamma$ is small enough. Note that by (\ref{smallness}), it holds $||\tilde{w}||_{H^{\frac{1}{2}}} \leq \epsilon(\gamma)$ with $\epsilon(\gamma) \longrightarrow 0$ as $\gamma \longrightarrow 0$. Let $F(x) = \frac{1}{2}x^2$ so that $\tilde{w}$ satisfies 
\begin{equation}
-\gamma \partial_x^{-2} \tilde{w} +H\tilde{w}_x  +\tilde{w} = F(\psi_{1,\gamma})- F(\tilde{\psi}_{1,\gamma})
\end{equation}
Denoting by $l_{1,\gamma}$ the operator defined in Proposition \ref{linoperator} with $c=1$, it holds:$$
l_{1,\gamma}= -\gamma \partial_x^{-2} + H\partial_x + 1 -F^{\prime}(\psi_{1,\gamma})
$$
and thus 
$$
l_{1,\gamma}\tilde{w}= F(\psi_{1,\gamma}) - F(\tilde{\psi}_{1,\gamma}) - F^{\prime}(\psi_{1,\gamma})(\psi_{1,\gamma} - \tilde{\psi}_{1,\gamma}) = \frac {\tilde{w}^2}{2}.
$$
\\
%with, for all $x\in \mathbb{R}$, $W_{\gamma}(x) \in [\psi_{1,\gamma}(x),\tilde{\psi}_{1,\gamma}(x)]$.\\
We will proceed by contradiction, suppose that $\tilde{w} \neq 0 $ and let $w = \frac{\tilde{w}}{||\tilde{w}||_{H^{\frac{1}{2}}}}$, we get 
\begin{align}
<l_{1,\gamma} w,w>_{L^2} &= \frac{||\tilde{w}||^3_{L^3}}{||\tilde{w}||^2_{H^{\frac{1}{2}}}}\nonumber\\
&\leq ||w||_{L^2}||\tilde{w}||_{L^2}~(\text{by} (\ref{L3L2}))\nonumber
\end{align}
so that
\begin{align}\label{smallnessoperator}
|<l_{1,0}w,w>_{L^2}|&= \left|<l_{1,\gamma}w,w>_{L^2} + <w^2,\psi_{1,\gamma}-Q_{1}>_{L^2} - \gamma||\partial_x^{-1}w||_{L^2}^2\right|\nonumber\\
&\lesssim \epsilon(\gamma) + \gamma.
\end{align}
Now, since $\frac{d}{dc} ||Q_{c}||^2_{L^2} >0$ then there exists $\alpha>0$ such that (see \cite{Bennett}): 
$$
<l_{c,0}u,u> \geq \alpha ||u||_{H^{\frac{1}{2}}}^2
$$
for all $u \in H^{\frac{1}{2}}(\mathbb{R})$ that satisfies the orthogonality conditions:\\
 $$<u,Q^{\prime}_{c}> = <u,Q_{c}>=0.$$ $w$ is even and $Q^{\prime}_{c}$ is odd, then we have the following orthogonality condition:
$$
<w,Q^{\prime}_{c}>_{L^2} = 0.
$$
Moreover, by differentiating
the equation (\ref{solitaryequationRBO}) with respect to c, and taking $c = 1$, we find that $l_{1,0}\partial_{c}|_{c=1}Q_{c}=-Q_{1}$ we get
\begin{align}
\left|<w,Q_{1}>\right| &= \left|<w,l_{1,0}\partial_{c}|_{c=1}Q_{c}>_{L^2}\right|\nonumber\\
&=\left|<l_{1,0}w,\partial_{c}|_{c=1}Q_{c}>_{L^2}\right|\nonumber\\
&\leq \left|<l_{1,\gamma}w,\partial_{c}|_{c=1}Q_{c}>_{L^2}\right| + \left|<w\partial_{c}|_{c=1}Q_{c},\psi_{c,\gamma}-Q_{1}>_{L^2}\right|\nonumber\\
&+\gamma\left|<w,\partial_x^{-2}\partial_{c}|_{c=1}Q_{c}>_{L^2}\right|\nonumber\\
&\lesssim (\epsilon(\gamma) + \gamma)||w||_{L^2}.
\end{align}
%where we have used that $\partial_{c}|_{c=1}Q_{c,0} \in L^{\infty}(\mathbb{R})$ in the last step. 
Then for $\gamma$ small enough, $w$ is almost orthogonal to $Q_{1}$. Therefore, we deduce that for $\gamma >0$ small enough,
$$
<l_{1,0}w,w> \geq \frac{\alpha}{2}||w||_{H^{\frac{1}{2}}} = \frac{\alpha}{2},
$$
which contradicts (\ref{smallnessoperator}) and proves the uniqueness result for $0<\gamma<\delta$. \\

Now let us prove the orbital stability result using the following proposition (see \cite{Bouard}):
\begin{proposition}\label{linoperator}
Let $\psi_{c,\gamma}$ be a solution of (\ref{solitaryequationRBO}) traveling with the speed c, and $l_{c,\gamma}$ be the linearized operator associated to the second derivative of the action functional E+cV at $\psi_{c,\gamma}$ defined by $l_{c,\gamma}v= -\gamma \partial_x^{-2}v +\mathcal{H}v_x + cv -\psi_{c,\gamma} v$. Assume that there exists $\alpha >0$ such that:$$
<l_{c,\gamma}v,v>_{L^2} \geq \alpha ||v||^2_{H^{\frac{1}{2}}(\mathbb{R})},
$$
for all $v \in H^{\frac{1}{2}}(\mathbb{R})$ satisfying the orthogonolaties
\begin{equation}
<v,\psi_{c,\gamma}>_{L^2}=<v,\psi_{c,\gamma}^{\prime}>_{L^2} = 0.
\end{equation}
Then $\psi$ is stable in $H^{\frac{1}{2}}(\mathbb{R})$.
\end{proposition}
Note that by remark (\ref{dilatation}) the orbital stability of $\psi_{1,\gamma}$ is equivalent to the orbital stability of $\psi_{\frac{\gamma_0}{\gamma},\gamma_0}$. To prove the orbital stability of this ground state for $\gamma >0$ small enough, we will proceed as following: first by the $H^{\frac{1}{2}}$ convergence result, it is easy to see that for $\gamma >0$ small enough, any $v\in H^{\frac{1}{2}}$ satisfying $<v,\psi_{1,\gamma}>=<v,\psi^{\prime}_{1,\gamma}>=0$ is almost orthogonal in $L^2(\mathbb{R})$ to $Q_{1}$ and $Q^{\prime}_{1}$, second $l_{1,0}$ is coercive in $H^{\frac{1}{2}}$ under these almost orthogonality conditions since  $ \frac{d}{dc}|_{c=1}||Q_{c}||^2_{L^2} >0$ (see \cite{Bennett}). Moreover:
\begin{align}
<l_{1,\gamma}v,v> &= \gamma ||\partial_x^{-1}v||_{L^2}^2 + <l_{1,0}v,v>_{L^2} - <v^2,\psi_{1,\gamma}-Q_{1}>_{L^2}\nonumber\\
&\geq \gamma||\partial_x^{-1}v||_{L^2}^2 + \alpha||v||^{2}_{H^{\frac{1}{2}}} - \epsilon(\gamma)||v||^2_{L^2}
\nonumber\\
&\geq \alpha||v||^2_{H^{\frac{1}{2}}}- \epsilon(\gamma)||v||^2_{H^{\frac{1}{2}}}\nonumber \\
&\geq \left(\alpha -\epsilon(\gamma)\right)||v||^2_{H^{\frac{1}{2}}}\nonumber
\end{align}
with $\alpha >0$ and $\epsilon(y) \longrightarrow 0$ as $y \longrightarrow 0$.\\
By proposition (\ref{linoperator}) we infer that there exists $\delta >0$ such that for all $\gamma \in ]0,\delta[$, $\psi_{1,\gamma}$ is orbitally stable.

\end{document}